\documentclass{article}
\usepackage{amssymb,latexsym,amsmath,amsthm,mathrsfs}
\usepackage{graphicx}
\newcommand{\comment}[1]{}

\newcommand{\Int}{{\textstyle \int}}
\begin{document}
\title{An observation on the sums
of divisors\footnote{Presented to the St. Petersburg Academy
on April 6, 1752. Originally published as
{\em Observatio de summis divisorum}, Novi
Commentarii Academiae scientiarum Imperialis Petropolitanae
\textbf{5} (1760), 59--74.
E243 in the Enestr{\"o}m index.
Translated from the Latin by Jordan Bell,
Department of Mathematics, University of Toronto, Toronto, Ontario, Canada.
Email: jordan.bell@gmail.com}}
\author{Leonhard Euler}
\date{}
\maketitle

\begin{center}
{\Large Summarium}
\end{center}
{\small
The divisors of a number are the numbers by which it can be divided without
a remainder; and among the divisors of any number are first unity and then
the number itself, since every number can be divided at least by unity and 
itself.
Now, those numbers which
besides unity and themselves have no other divisors we call prime,
while the others, which permit division without remainder
by some number other than themselves, we call composite;
and in common Arithmetic a method is taught for
finding all the divisors of any number.
The Author in this dissertation considers 
the sum of all divisors of any given number, not to the ends that others
typically pursue,
for the investigation of perfect or amicable numbers
and other questions of this kind,
but to explore the order and as it were the law by which the sums
of the divisors of each number proceed conveniently.
This should certainly be seen as most concealed, since for prime numbers
the sum of divisors will exceed them by unity,
but for composites it is the greater the
more prime factors that comprise them.
Therefore, since the rule for the progression of the prime numbers
is so far a great mystery, into which not even Fermat was able to penetrate,
since the rule for these is clearly involved in the sums of divisors,
who should doubt that these too are not subject to any law?
So this dissertation merits all the more attention, because such a law
is brought here to light, even if it has not yet been demonstrated with complete
rigor.
The same thing happens for the Author as before with the theorem of Fermat,
that soon after the defect in the demonstration has been corrected.\footnote{Translator: Right now I presume this refers to Euler's E228 and E241, both
on representing integers as sums of two square. However I have
not read these papers.} 
For what is still desired in the demonstration given here
will be at once supplied in the following dissertation.

In order to clearly explain this, the Author uses the sign $\Int$
to indicate the sum of the divisors of any number. Thus
$\Int n$ indicates the sum of all the divisors of the number $n$.
One can see that:
\[
\begin{array}{lll}
\Int 1=1,&&\Int 11=1+11=12,\\
\Int 2=1+2=3,&&\Int 12=1+2+3+4+6+12=28,\\
\Int 3=1+3=4,&&\Int 13=1+13=14,\\
\Int 4=1+2+4=7,&&\Int 14=1+2+7+14=24,\\
\Int 5=1+5=6,&&\Int 15=1+3+5+15=24,\\
\Int 6=1+2+3+6=12,&&\Int 16=1+2+4+8+16=31,\\
\Int 7=1+7=8,&&\Int 17=1+17=18,\\
\Int 8=1+2+4+8=15,&&\Int 18=1+2+3+6+9+18=39,\\
\Int 9=1+3+9=13,&&\Int 19=1+19=20,\\
\Int 10=1+2+5+10=18,&&\Int 20=1+2+4+5+10+20=42\\
&\textrm{etc.}&
\end{array}
\]
These sums can be easily defined from the
known principle\footnote{Translator: cf. pp. 53--54 of
Andr\'e
Weil,
{\em Number theory: an approach through history}.} that
the sum of the divisors of the product $mnpq$, where the factors $m,n,p,q$
are taken to be prime to each other, is equal to the product of each of
the sums, or
\[
\Int mnpq=\Int m\cdot \Int n\cdot \Int p\cdot \Int q.
\]
Thus
\[
\Int 20=\Int 4\cdot 5=\Int 4\cdot \Int 5=7\cdot 6=42 \quad
\textrm{and} \quad
\Int 360=\Int 8\cdot 9\cdot 5=\Int 8\cdot \Int 9\cdot \Int 5=
15\cdot 13\cdot 6=1170.
\]
The Author also considers the correspondence between the sums of divisors
and the numbers in their usual order, which proceed as follows
\[
\begin{array}{llllllllllllllllll}
\textrm{number}&1,&2,&3,&4,&5,&6,&7,&8,&9,&10,&11,&12,&13,&14,&15,&16,&\textrm{etc.},\\
\textrm{sum of div.}&1,&3,&4,&7,&6,&12,&8,&15,&13,&18,&12,&28,&14,&24,&24,&31,&\textrm{etc.},
\end{array}
\]
in which progression certainly no law is seen, since now it is greater
and now it is less, now even and now odd, and especially, the order of
the prime numbers is clearly involved; since this is so impenetrable, who
would suspect a law in this series? But nevertheless, the Author
shows that these numbers constitute a series of the kind that is usually
called recurrent,
such that each term can be determined from some of the preceding
according to a certain law.
As $\Int n$ denotes the sum of the divisors of the number $n$,
the expression $\Int (n-a)$ will denote the sum of the divisors
of the number $n-a$. With this notation established, the law
found by the Author is that
\[
\begin{split}
&\Int n=\Int (n-1)+\Int (n-2)-\Int (n-5)-\Int (n-7)+\Int (n-12)
+\Int (n-15)\\
&-\Int (n-22)-\Int (n-26)+\Int (n-35)+\Int (n-40)-\textrm{etc.},
\end{split}
\]
in which rule the signs are taken two $+$ followed by two $-$; and
the numbers $1,2,5,7,12,15$ which are continually subtracted from $n$
can be easily found from their differences:
\[
\begin{array}{ll}
\textrm{numbers}&1,2,5,7,12,15,22,26,35,40,51,57 \quad \textrm{etc.},\\
\textrm{differences}&1,3,2,5,3,7,4,9,5,11,6\quad  \textrm{etc.},
\end{array}
\]
which come alternately from two sequences.
The relation is conveniently represented in the following way:
\[ \Int n = \begin{cases}
+\Int (n-1)-\Int(n-5)+\Int(n-12)-\Int(n-22)+\Int(n-35)-\textrm{etc.}&\\
+\Int (n-2)-\Int(n-7)+\Int(n-15)-\Int(n-26)+\Int(n-40)-\textrm{etc.}&
\end{cases}
\]
For the application of this formula to any number, one should understand
that these terms are taken until the number after the $\Int$ sign is negative,
which is omitted together with all further terms. As well then, if
the term $\Int (n-n)$ or $\Int 0$ occurs, since this by itself is not
determinate, in this case the number $n$ itself is written in place of the
term. Thus, according to this law it will be
\[
\Int 21=\Int 20+\Int 19 - \Int 16-\Int 14+\Int 9+\Int 6
\]
or
\[
\Int 21=42+20-31-24+13+12=87-55=32;
\]
and
\[
\Int 22=\Int 21+\Int 20-\Int 17-\Int 15+\Int 10+\Int 7-\Int 0
\]
or
\[
\Int 22=32+42-18-24+18+8-22=100-64=36.
\]
The Author has deduced the marvelous law for this progression from the consideration
of the product
\[
(1-x)(1-x^2)(1-x^3)(1-x^4)(1-x^5)(1-x^6)(1-x^7)\, \textrm{etc.}
\]
whose factors continue to infinity. But if this is expanded by actual
multiplication it will be found to form the following series
\[
1-x-x^2+x^5+x^7-x^{12}-x^{15}+x^{22}+x^{26}-x^{35}-x^{40}+\textrm{etc.},
\]
and it indeed has been possible to continue this operation further,
whence the law of this series and the progression of the exponents has
been concluded at least by induction, which seems perhaps to suffice for many.
Truly the Author candidly admits that this observed agreement to be
by no means demonstrated, but that its demonstration is still desired,
which however he has communicated hardly much later to the Academy.
With the equality of the product and the expanded series, the mentioned
theorem about the order in the sums of divisors is thence clearly
demonstrated, such that there cannot be any more doubt, even if
it is accomplished by logarithms and differentiation, which seem to have
little to do with the nature of divisors. From this, one may
see how closely and wonderfully infinitesimal Analysis is connected not
only with usual Analysis, but even with the theory of numbers, which
seems remote from that higher kind of calculus.}

1. For a given number $n$, let the formula $\Int n$ denote the sum of all
the divisors of the number $n$. Then since unity has no other divisors
besides itself, $\Int 1=1$; and since a prime number has exactly
two factors, unity and itself, if $n$ is a prime number then
$\Int n=1+n$. Next, since a perfect number is equal to the sum
of its aliquot parts, where aliquot parts are the divisors other than
itself, it is immediate that the sum of the divisors of a perfect number
will be twice the number itself;
that is, if $n$ is a perfect number then $\Int n=2n$.
Further, usually a number is called abundant if the sum of its
aliquot parts is greater than it, so if $n$ is an abundant
number, then $\Int n>2n$; and if $n$ is a deficient number, that is such
that the sum of its aliquot parts is less than it, then $\Int n<2n$.

2. In this manner now, one can easily express with signs the
essences of numbers, as far as this is contained in the sum of
the aliquot parts or divisors. For if $\Int n=1+n$, then $n$ will
be a prime number, if $\Int n=2n$, then $n$ will be a perfect number,
and if $\Int n>2n$ or $\Int n<2n$, then $n$ will be an abundant
or deficient number respectively. Here the question can also be dealt with
of those
numbers which are usually called amicable, that is, for which the sum of
the aliquot parts of each is equal to the other. For if $m$ and
$n$ are amicable numbers, since
the sum of the aliquot parts of the number $m$ is $=\Int m-m$
and of the number $n$ is $=\Int n-n$, by the nature of these numbers
it will be $n=\Int m-m$ and $m=\Int n-n$, and thus $\Int m=\Int n=m+n$.
Hence two amicable numbers have the same sum of divisors, which is
simultaneously equal to the sum of both the numbers.

3. To easily find the sum of the divisors of any given number, it is most
convenient to resolve the number into two factors which are
prime to each other. For if $p$ and $q$ are numbers that are prime
to each other, or which besides unity have no common divisor, then the sum
of the divisors of the product $pq$ will be equal to the product
of the sums of the divisors of both, or
\[
\Int pq=\Int p\cdot \Int q.
\]
Thence, with the sums of the divisors of smaller numbers found, it will
not be difficult to extend the discovery of the sum of divisors
to greater numbers.

4. If $a,b,c,d$ etc. are the prime numbers, then every number, no matter
what size, can always be reduced to the form $a^\alpha b^\beta c^\gamma d^\delta$ etc.;
having gotten this form, the sum of the divisors of this number or
$\Int a^\alpha b^\beta c^\gamma d^\delta$ etc. will be
\[
=\Int a^\alpha \cdot \Int b^\beta \cdot \Int c^\gamma \cdot \Int d^\delta\cdot
\, \textrm{etc.}
\]
But because $a,b,c,d$ etc. are prime numbers,
\[
\Int a^\alpha = 1+a+a^2+\cdots+a^\alpha=\frac{a^{\alpha+1}-1}{a-1}
\]
and so
\[
\Int a^\alpha b^\beta c^\gamma d^\delta \, \textrm{etc.}=\frac{a^{\alpha+1}-1}{a-1}
\cdot \frac{b^{\beta+1}-1}{b-1}\cdot \frac{c^{\gamma+1}-1}{c-1}
\cdot \frac{d^{\delta+1}-1}{d-1}\cdot \, \textrm{etc.}
\]
It will therefore suffice to have found only the sums of the divisors
of all the powers of prime numbers.\footnote{Translator: The {\em Opera omnia}
cites several sources on the multiplicativity of the sum of divisors. 
In particular, it mentions Euler's paper E152, on amicable numbers.}

5. However, I shall not pursue this inquiry further. To come nearer to
what I have set to treat, let me write out for later use the sums of the divisors
of numbers proceeding according to the natural order.

\begin{tabular}{p{3cm}p{3cm}p{3cm}p{3cm}}
${\int}1=1$&${\int}26=42$&${\int}51=72$&${\int}76=140$\\
${\int}2=3$&${\int}27=40$&${\int}52=98$&${\int}77=96$\\
${\int}3=4$&${\int}28=56$&${\int}53=54$&${\int}78=168$\\
${\int}4=7$&${\int}29=30$&${\int}54=120$&${\int}79=80$\\
${\int}5=6$&${\int}30=72$&${\int}55=72$&${\int}80=186$\\
${\int}6=12$&${\int}31=32$&${\int}56=120$&${\int}81=121$\\
${\int}7=8$&${\int}32=63$&${\int}57=80$&${\int}82=126$\\
${\int}8=15$&${\int}33=48$&${\int}58=90$&${\int}83=84$\\
${\int}9=13$&${\int}34=54$&${\int}59=60$&${\int}84=224$\\
${\int}10=18$&${\int}35=48$&${\int}60=168$&${\int}85=108$\\
${\int}11=12$&${\int}36=91$&${\int}61=62$&${\int}86=132$\\
${\int}12=28$&${\int}37=38$&${\int}62=96$&${\int}87=120$\\
${\int}13=14$&${\int}38=60$&${\int}63=104$&${\int}88=180$\\
${\int}14=24$&${\int}39=56$&${\int}64=127$&${\int}89=90$\\
${\int}15=24$&${\int}40=90$&${\int}65=84$&${\int}90=234$\\
${\int}16=31$&${\int}41=42$&${\int}66=144$&${\int}91=112$\\
${\int}17=18$&${\int}42=96$&${\int}67=68$&${\int}92=168$\\
${\int}18=39$&${\int}43=44$&${\int}68=126$&${\int}93=128$\\
${\int}19=20$&${\int}44=84$&${\int}69=96$&${\int}94=144$\\
${\int}20=42$&${\int}45=78$&${\int}70=144$&${\int}95=120$\\
${\int}21=32$&${\int}46=72$&${\int}71=72$&${\int}96=252$\\
${\int}22=36$&${\int}47=48$&${\int}72=195$&${\int}97=98$\\
${\int}23=24$&${\int}48=124$&${\int}73=74$&${\int}98=171$\\
${\int}24=60$&${\int}49=57$&${\int}74=114$&${\int}99=156$\\
${\int}25=31$&${\int}50=93$&${\int}75=124$&${\int}100=217$
\end{tabular}

6. If we contemplate now the series of the numbers 1, 3, 4, 7, 6, 12, 8, 15,
13, 18, 12, 28 etc.
which the sums of divisors corresponding to the numbers proceeding in their
natural order constitute, not only is there no apparent law for the progression,
but the order of these numbers seems so disturbed that they seem
to be bound by no law whatsoever.
For this series is clearly mixed up with the order of prime numbers,
since the term of index $n$, or $\Int n$, will always
be $=n+1$ exactly when $n$ is a prime number;
but it is well known that thus far it has not been possible to 
refer the prime numbers to any certain law of progression.  
And since our series involves the rule not only of the prime numbers
but also all the other numbers, insofar as they are composed from primes,
its law would seem more difficult to find than that of the series
of prime numbers alone.

7. Since this is the case, I seem to have advanced the science of numbers
by not a small amount when I found a certain fixed law according
to which the terms of the given series $1,3,4,7,6$ etc. progress,
such that by this law each term of the series can be defined from the preceding;
for I have found, which seems rather wonderful, that this series belongs
to the kind of progression which are usually called recurrent and
whose nature is such that each term is determined from the preceding
according to some certain rule of relation. And who would have ever believed
that this series which is so disturbed and which seems to have nothing
in common with recurrent series would nevertheless be included
in this type of series, and that it would be possible to assign a
scale of relation for it?\footnote{Translator: scalam relationis=scale of relation=recurrence relation. Euler uses this term in his November 10, 1742 letter to Nicolaus
I Bernoulli.}

8. Since the term of this series corresponding to the index $n$, which
indicates the sum of the divisors of the number $n$, is $=\Int n$, the
prior terms in descending order are $\Int (n-1), \Int (n-2), \Int (n-3),
\Int (n-4), \Int (n-5)$ etc. And any term of this series, namely
$\Int n$, is conflated from some of the prior terms, as
\[
\begin{split}
&\Int n=\Int (n-1)+\Int (n-2)-\Int(n-5)-\Int(n-7)+\Int(n-12)+\Int(n-15)\\
&-\Int(n-22)-\Int(n-26)+\Int(n-35)+\Int(n-40)-\Int(n-51)-\Int(n-57)\\
&+\Int(n-70)+\Int(n-77)-\Int(n-92)-\Int(n-100)+\Int(n-117)+\Int(n-126)-\textrm{etc.}
\end{split}
\]
Or since the signs $+$ and $-$ occur alternately in pairs,
this series can be separated easily into two like this:
\[
\Int n =\begin{cases}
\Int(n-1)-\Int(n-5)+\Int(n-12)-\Int(n-22)+\Int(n-35)-\Int(n-51)+\textrm{etc.}&\\
\Int(n-2)-\Int(n-7)+\Int(n-15)-\Int(n-26)+\Int(n-40)-\Int(n-57)+\textrm{etc.}&
\end{cases}
\]

9. From the above form, the order of the numbers which are successively
subtracted from $n$ in each series is easily seen;
for each series is of the second order, having constant
second differences.\footnote{Translator: Namely a second order
arithmetic progression.}
In fact, the numbers of the first series together with both their
first and second differences are

{\scriptsize
\[
\begin{array}{lllllllllllllllllll}
&1,&&5,&&12,&&22,&&35,&&51,&&70,&&92,&&117,&\textrm{etc.},\\
\textrm{1st diff.}&&4,&&7,&&10,&&13,&&16,&&19,&&22,&&25,&&\textrm{etc.},\\
\textrm{2nd diff.}&&&3,&&3,&&3,&&3,&&3,&&3,&&3&&\textrm{etc.}&
\end{array}
\]
}
Whence the general term of this series\footnote{Translator: Let $\Delta f(x)=f(x+1)-f(x)$ and $\Delta^{k+1} f(x)=
\Delta^k f(x+1)-\Delta^k f(x)$. If $f(x)$ is a polynomial, then
\[
f(x+a)=\sum_{k=0}^\infty \binom{x}{k} \Delta^k f(a).
\]
This is sometimes called ``Newton's series'' for $f$.}
is $=\frac{3xx-x}{2}$ and thus contains exactly the pentagonal numbers.
The other series is
{\scriptsize
\[
\begin{array}{lllllllllllllllllll}
&2,&&7,&&15,&&26,&&40,&&57,&&77,&&100,&&126,&\textrm{etc.},\\
\textrm{1st diff.}&&5,&&8,&&11,&&14,&&17,&&20,&&23,&&26,&&\textrm{etc.},\\
\textrm{2nd diff.}&&&3,&&3,&&3,&&3,&&3,&&3,&&3&&\textrm{etc.}&
\end{array}
\]
}
and hence the general term is $\frac{3xx+x}{2}$ and contains the series of
pentagonal numbers continued backwards.\footnote{Translator: Continued to
negative indices? viz. $\frac{-x(-3x-1)}{2}=\frac{3x+1}{2}$.}

10. It is highly noteworthy here that the series of pentagonal numbers, itself
and continued backwards, is applied with great effect to the order of
the series of sums of divisors, since certainly one
would not at all suspect there to be a connection 
between the pentagonal numbers and sums of divisors.
For if one writes the series of pentagonal numbers forwards
and continued backwards in this way
\[
\textrm{etc.} \quad 77, 57, 40, 26, 15, 7, 2, 0, 1, 5, 12, 22, 35, 51, 70, 92
\quad \textrm{etc.}
\]
our formula enclosing the order of the sums of divisors can be presented
with alternating signs ordered in this way
\[
\begin{split}
&\textrm{etc.}-\Int (n-15)+\Int(n-7)-\Int(n-2)+\Int( n-0)-\Int(n-1)\\
&+\Int(n-5)-\Int(n-12)+\Int(n-22)-\textrm{etc.}=0,
\end{split}
\]
in which the series on both sides continue to infinity, but in every case if
it is correctly applied to our use a determinate numbers of terms will arise.

11. For if we want to find the sum of the divisors of the number $n$
by means of our first exhibited formula
\[
\begin{split}
&\Int n=\Int(n-1)+\Int(n-2)-\Int(n-5)-\Int(n-7)+\Int(n-12)+\Int(n-15)\\
&-\Int(n-22)-\Int(n-26)+\Int(n-35)+\Int(n-40)-\Int(n-51)-\Int(n-57)\\
&+\Int(n-70)+\Int(n-77)-\Int(n-92)-\Int(n-100)+\textrm{etc.},
\end{split}
\]
with the sums of the divisors of smaller numbers known, then we only need
to take the terms in this formula until we reach sums of the divisors
of negative numbers. Namely, all the terms which contain negative numbers
after the $\Int$ sign are rejected; whence it is clear that if $n$ is a small
number just a few terms suffice, while if $n$ is a larger number then it will
be necessary to take more terms from our general formula.

12. Therefore, the sum of the divisors of a given number $n$ is composed
from the sum of divisors of some smaller numbers which I assume to be
known, since in each case the sums for negative numbers are rejected.
This is an easy provision, because one cannot even take the sum of the divisors
of negative numbers; but it should be explained how this operation
is done in those cases in which our formula yields the term $\Int(n-n)$
or $\Int 0$, which, since zero is divisible by all numbers, seems either
infinite or indeterminate. This case will occur exactly when $n$ is
a number from either the series of pentagonal numbers
or the series continued backwards; then in these cases the number $n$
itself should be taken in place of the term $\Int(n-n)$ or $\Int 0$,
and should be written with the sign which the term $\Int(n-n)$ is
affixed with in our formula.

13. With these precepts for the use of our formula explained,
to begin with I shall give examples with small numbers which can be
easily examined by means of our formula, and
simultaneously the truth of the formula will be recognized.

\begin{align*}
&\Int 1=\Int 0\\
\textrm{or} \quad &\Int 1=1=1\\
&\Int 2=\Int 1+\Int 0\\
\textrm{or} \quad &\Int 2=1+2=3\\
&\Int 3=\Int 2+\Int 1\\
\textrm{or} \quad &\Int 3=3+1=4\\
&\Int 4=\Int 3+\Int 2\\
\textrm{or} \quad &\Int 4=4+3=7\\
&\Int 5=\Int 4+\Int 3-\Int 0\\
\textrm{or} \quad &\Int 5=7+4-5=6\\
&\Int 6=\Int 5+\Int 4-\Int 1\\
\textrm{or} \quad &\Int 6=6+7-1=12\\
&\Int 7=\Int 6+\Int 5-\Int 2 - \Int 0\\
\textrm{or} \quad &\Int 7=12+6-3-7=8\\
&\Int 8=\Int 7+\Int 6-\Int 3-\Int 1\\
\textrm{or} \quad &\Int 8=8+12-4-1=15\\
&\Int 9=\Int 8+\Int 7-\Int 4-\Int 2\\
\textrm{or} \quad &\Int 9=15+8-7-3=13\\
&\Int 10=\Int 9+\Int 8-\Int 5-\Int 3\\
\textrm{or} \quad &\Int 10=13+15-6-4=18\\
&\Int 11=\Int 10+\Int 9-\Int 6-\Int 4\\
\textrm{or} \quad &\Int 11=18+13-12-7=12\\
&\Int 12=\Int 11+\Int 10-\Int 7-\Int 5+\Int 0\\
\textrm{or} \quad &\Int 12=12+18-8-6+12=28.
\end{align*}

14. By inspecting these examples with attention and also by continuing
to greater numbers, it will be apparent not without admiration how, as
it were against expectation, that the true sum of divisors of the given
number is obtained;
and to make it easier to recognize this pattern, I have already given
above the sums of the divisors of all numbers not greater than one hundred,
whence the truth of our formula can be tested with greater numbers.
In particular we will find not without delight that the given number is
prime when the sum found from our formula for it is greater than the number
by unity. Let us work out an example to this end, with the given
number $n=101$, and test it as if ignorant about whether or not
this number is prime. The operation will happen thus:
\[
\begin{array}{rcrrrrrrrr}
\Int 101&=&\Int 100&+\Int 99&-\Int 96&-\Int 94&+\Int 89&+\Int 86&-\Int 79&-\Int 75\\
&&217&+156&-252&-144&+90&+132&-80&-124\\
&&+\Int 66&+\Int 61&-\Int 50&-\Int 44&+\Int 31&+\Int 24&-\Int 9&-\Int 1\\
&&+144&+62&-93&-84&+32&+60&-13&-1.
\end{array}
\]
Therefore by collecting the pairs of two terms together we will have
\begin{eqnarray*}
\Int 101&=&+373-396\\
&&+222-204\\
&&+206-177\\
&&+92-14
\end{eqnarray*}
or
\[
\Int 101=+893-791=102.
\]
Therefore, we find that the sum of the divisors of the number 101 is greater
than it by unity, namely 102,
whence even if it were not otherwise known, it clearly follows that
the number 101 is prime. This rightly seems miraculous, since no operation
was done which referred in any way to the calculation of divisors;
also, the divisors whose sum is found by this method remain themselves
unknown, although they can frequently be figured out from the consideration
of this sum.\footnote{Translator: Is Euler saying that if we know
$\Int n$ for all $n$, then for any particular number $n$ it is simple
to find its divisors just using operations involving $\Int$?}

15. These special properties which the sums of divisors are gifted with
would be no less memorable if their demonstration were obvious, and as it
were exposed to the daylight. But the demonstration was in fact abstruse
and depended on rather difficult properties of numbers,
whence to no small degree the value of this law discovered for the progression
is increased; for the investigation of truths is to be recommended the more
the more hidden they are. Truly, I am compelled to admit that now not only
have I not been able to find a demonstration of this truth, but that I have
even nearly been brought to despair, and I do not know whether because
of this the knowledge of a truth
whose demonstration is hidden to us
should be valued even more highly.
And so this truth has been confirmed by a great many examples, since it has
not been permitted that I exhibit a demonstration of it.

16. Thus here we have an extraordinary example of the kind of proposition
whose truth we can in no way doubt, even if we have not achieved its
demonstration. This will seem rather surprising to most, since
in common mathematics no propositions are counted as true
unless they can be derived from indubitable principles. Yet in the meanwhile,
I have come to the knowledge of this truth not by chance and, as it
were, by divination; for to whom would it have come to mind to try to
elicit by conjecture alone an order that might perhaps occur in the sums of divisors, from
the nature of recurrent series and of the pentagonal numbers?
For which reason I judge it not to be foreign from our purpose
if I clearly explain the way by which came to the knowledge of
this order, especially since it is very recondite and was discovered
in a long roundabout way.

17. I was led to this observation by considering the infinite formula
\[
s=(1-x)(1-x^2)(1-x^3)(1-x^4)(1-x^5)(1-x^6)(1-x^7)(1-x^8) \, \textrm{etc.},
\]
which if actually expanded by multiplication and then arranged according
to the powers of $x$, I discovered to be transformed into the following series
\[
s=1-x-x^2+x^5+x^7-x^{12}-x^{15}+x^{22}+x^{26}-x^{35}-x^{40}+x^{51}+x^{57}
-\textrm{etc.},
\]
where exactly those numbers occur in the exponents of $x$ which
I described above, namely the pentagonal numbers themselves and them
continued backwards. To more easily see this order, the series can be
exhibited thus, going to infinity on each side
\[
s=\textrm{etc.}+x^{26}-x^{15}+x^{7}-x^{2}+x^0-x^1+x^5-x^{12}+x^{22}-x^{35}
+x^{51}-\textrm{etc.}
\]

18. The equality of these two formulas exhibiting $s$ is now the very
thing which I am not able to confirm with a solid demonstration;
nevertheless
undertaking to successively multiply out the factors of the first formula
\[
s=(1-x)(1-x^2)(1-x^3)(1-x^4)(1-x^5) \, \textrm{etc.}
\]
leads to the initial terms of the other series
\[
s=1-x-x^2+x^5+x^7-x^{12}-x^{15}+\textrm{etc.},
\]
and neither is it difficult to see that the two signs $+$ and $-$
occur alternately in pairs and that the exponents of the powers of
$x$ follow the same law that I explained enough already.
But conceded the equality of these two infinite formulas, the properties
of the sums of divisors which I indicated before can be rigidly
demonstrated; and on the other hand, if these properties are admitted as
true, the true agreement of our two formulas will follow.

19. If for doing the demonstration we assume that both
\[
s=(1-x)(1-x^2)(1-x^3)(1-x^4)(1-x^5) \, \textrm{etc.}
\]
and
\[
s=1-x-x^2+x^5+x^7-x^{12}-x^{15}+x^{22}+x^{26}-\textrm{etc.},
\]
by taking logarithms we get
\[
ls=l(1-x)+l(1-x^2)+l(1-x^3)+l(1-x^4)+l(1-x^5)+\textrm{etc.}
\]
and
\[
ls=l(1-x-x^2+x^5+x^7-x^{12}-x^{15}+x^{22}+x^{26}-\textrm{etc.}).
\]
Then taking the differentials of each formula we have
\[
\frac{ds}{s}=-\frac{dx}{1-x}-\frac{2xdx}{1-x^2}-\frac{3x^2dx}{1-x^3}
-\frac{4x^3dx}{1-x^4}-\frac{5x^4dx}{1-x^5}-\textrm{etc.}
\]
and
\[
\frac{ds}{s}=\frac{-dx-2xdx+5x^4dx+7x^6dx-12x^{11}dx-15x^{14}dx+\textrm{etc.}}{1-x-x^2+x^5+x^7-x^{12}-x^{15}+x^{22}+x^{26}-\textrm{etc.}}
\]
Let us multiply both of these by $\frac{-x}{dx}$, so that we have
\begin{align*}
\textrm{I.}&-\frac{xds}{dx}=\frac{x}{1-x}+\frac{2x^2}{1-x^2}+\frac{3x^3}{1-x^3}
+\frac{4x^4}{1-x^4}+\frac{5x^5}{1-x^5}+\textrm{etc.}\\
\textrm{II.}&-\frac{xds}{sdx}=\frac{x+2x^2-5x^5-7x^7+12x^{12}+15x^{15}
-22x^{22}-26x^{26}+\textrm{etc.}}{1-x-x^2+x^5+x^7-x^{12}-x^{15}+x^{22}+x^{26}
-\textrm{etc.}}
\end{align*}

20. First, let's consider the first of the two equal expressions,
and let us convert all the terms into geometric progressions in
usual manner; with this done, arranging these 
infinitely many geometric progressions according to powers of $x$ will yield:
\[
\begin{array}{rlllllllllllll}
-\frac{xds}{sdx}=&x^1&+x^2&+x^3&+x^4&+x^5&+x^6&+x^7&+x^8&+x^9&+x^{10}&+x^{11}&+x^{12}&+\textrm{etc.}\\
&&+2&&+2&&+2&&+2&&+2&&+2&\\
&&&+3&&&+3&&&+3&&&+3&\\
&&&&+4&&&&+4&&&&+4&\\
&&&&&+5&&&&&+5&&&\\
&&&&&&+6&&&&&&+6&\\
&&&&&&&+7&&&&&&\\
&&&&&&&&+8&&&&&\\
&&&&&&&&&+9&&&&\\
&&&&&&&&&&+10&&&\\
&&&&&&&&&&&+11&&\\
&&&&&&&&&&&&+12&
\end{array}
\]

21. Now if the coefficients of all the powers of $x$ are collected, one
will have
\[
-\frac{xds}{sdx}=x^1+x^2(1+2)+x^3(1+3)+x^4(1+2+4)+x^5(1+5)+x^6(1+2+3+6)+
\textrm{etc.},
\]
where it is clear that the coefficient of each power of $x$ is the sum of
all the numbers by which the exponent of the power is divisible.
Namely, the coefficient of the power $x^n$ will be the sum of all the divisors
of the number $n$;
thus according to the manner of signification explained above it will be
$=\Int n$. Then 
the series found equal to $-\frac{xds}{sdx}$ can thus be exhibited as
\[
-\frac{xds}{sdx}=x\Int 1+x^2\Int 2+x^3\Int 3+x^4\Int 4+x^5\Int 5+
x^6\Int 6+x^7\Int 7+\textrm{etc.},
\]
and by putting $x=1$ this yields the progression of the sums of divisors,
which assembles all all numbers proceeding in the natural order.

22. Let us now designate this series by $t$, so that
\[
t=x^1\Int 1+x^2\Int 2+x^3\Int 3+x^4\Int 4+x^5\Int 5+x^6\Int 6+x^7\Int 7+\textrm{etc.},
\]
and as $t=-\frac{xds}{sdx}$, it will also be that
\[
t=\frac{x^1+2x^2-5x^5-7x^7+12x^{12}+15x^{15}-22x^{22}-26x^{26}+\textrm{etc.}}{1-x-x^2+x^5+x^7-x^{12}-x^{15}+x^{22}+x^{26}-\textrm{etc.}}
\]
Then it is necessary that the series obtained for $t$ from the expansion
of this fraction be equal to that which the prior form has provided.
From this it is apparent that that series found for $t$ is recurrent:
each of its terms is determined from the preceding by a certain scale
of relation, which the denominator $1-x-x^2+x^5+x^7-\textrm{etc.}$
indicates.

23. Now so that the character of this recurrent series can be easily understood,
let us equate the two values found for $t$ and in order to get rid of
the fraction, let each be multiplied by the denominator $1-x-x^2+x^5+x^7-x^{12}
-x^{15}+\textrm{etc.}$ 
This done, by arranging the terms according to powers of $x$ this will
arise
{\scriptsize
\[
\setlength{\arraycolsep}{.5\arraycolsep}
\begin{array}{rrrrrrrrrrrrr}
x^1\Int 1&+x^2\Int 2&+x^3\Int 3&+x^4\Int 4&+x^5\Int 5&+x^6\Int 6&+x^7\Int 7
&+x^8\Int 8&+x^9\Int 9&+x^{10}\Int 10&+x^{11}\Int 11&+x^{12}\Int 12&+\textrm{etc.}\\
&-\Int 1&-\Int 2&-\Int 3&-\Int 4&-\Int 5&-\Int 6&-\Int 7&-\Int 8&-\Int 9&-\Int 10&-\Int 11&\\
&&-\Int 1&-\Int 2&-\Int 3&-\Int 4&-\Int 5&-\Int 6&-\Int 7&-\Int 8&-\Int 9&-\Int 10&\\
&&&&&+\Int 1&+\Int 2&+\Int 3&+\Int 4&+\Int 5&+\Int 6&+\Int 7&\\
&&&&&&&+\Int 1&+\Int 2&+\Int 3&+\Int 4&+\Int 5&\\
&&&&&&&&&&&\vdots&\\
=x^1&+2x^2&*&*&-5x^5&*&-7x^7&*&*&*&*&+12x^{12}&+\textrm{etc.}
\end{array}
\]
}

24. Since now the coefficients of each power of $x$ need to destroy
each other, we may elicit the following equalities
\begin{eqnarray*}
\Int 1=1,&&\Int 7=\Int 6+\Int 5-\Int 2-7,\\
\Int 2=\Int 1+2,&&\Int 8=\Int 7+\Int 6-\Int 3-\Int 1,\\
\Int 3=\Int 2+\Int 1,&&\Int 9=\Int 8+\Int 7-\Int 4-\Int 2,\\
\Int 4=\Int 3+\Int 2,&&\Int 10=\Int 9+\Int 8-\Int 5-\Int 3,\\
\Int 5=\Int 4+\Int 3-5,&&\Int 11=\Int 10+\Int 9-\Int 6-\Int 4,\\
\Int 6=\Int 5+\Int 4-\Int 1,&&\Int 12=\Int 11+\Int 10-\Int 7-\Int 5+12\\
&\textrm{etc.},&
\end{eqnarray*}
which clearly reduce to these
\begin{eqnarray*}
\Int 1=1,&&\Int 7=\Int (7-1)+\Int (7-2)-\Int (7-5)-7,\\
\Int 2=\Int (2-1)+2,&&\Int 8=\Int (8-1)+\Int (8-2)-\Int (8-5)-\Int (8-7),\\
\Int 3=\Int (3-1)+\Int (3-2),&&\Int 9=\Int (9-1)+\Int (9-2)-\Int (9-5)-\Int (9-7),\\
\Int 4=\Int (4-1)+\Int (4-2),&&\Int 10=\Int (10-1)+\Int (10-2)-\Int (10-5)-\Int (10-7),\\
\Int 5=\Int (5-1)+\Int (5-2)-5,&&\Int 11=\Int (11-1)+\Int (11-2)-\Int(11-5)-\Int(11-7),\\
\Int 6=\Int (6-1)+\Int(6-2)-\Int (6-5),&&\Int 12=\Int(12-1)+\Int(12-2)-\Int(12-5)-\Int(12-7)+12.
\end{eqnarray*}

25. Here it is evident that the numbers which ought to be continually
subtracted from the given number, the sum of whose divisors is sought,
are the very numbers from the series $1,2,5,7,12,15,22,26$ etc.;
in each case, they are to be taken as long as they do not exceed the given
number. As well the signs follow the same rule which was described above.
Therefore for any given number $n$ it will clearly be
\[
\Int n=\Int (n-1)+\Int (n-2)-\Int(n-5)-\Int(n-7)+\Int(n-12)+\Int(n-15)-\textrm{etc.};
\] 
the terms are to be continued until the numbers having the sign
$\Int$ in front of them become negative.
So from the origin of this recurrent series the rule is transparent
why this progression is not continued any further.

26. Then, for what pertains to the actual numbers which are appended
at the end of certain of the found formulas,
it is clear that they arise from the numerator of the fraction
whose value was found expressing $t$ (\S 22),
and interrupt the law of continuity for exactly those cases in which
the number $n$ is a term from the series $1,2,5,7,12,15,22,26$ etc.,
but even in this case the law of signs is not affected. 
In these cases the actual number which is to be added
is equal to the given number itself, keeping the same sign;
and if we consider the law described before, we see that this
number corresponds to the term $\Int(n-n)$ there. From this
the rule is transparent why, whenever in applying the formula
\[
\Int n=\Int (n-1)+\Int(n-2)-\Int(n-5)-\Int(n-7)+\Int(n-12)+\textrm{etc.}
\]
the term $\Int(n-n)$ is encountered, it is not omitted, but rather
its the number $n$ itself should be written for its value.
Therefore the rule explained above is confirmed in all
its parts.

\end{document}